\documentclass[12pt]{article}
\usepackage{amssymb}
\usepackage{amsmath}

\begin{document}

\begin{center}
{\bfseries On gaps between zeros of the Riemann zeta-function}
\end{center}
\begin{center}
Shaoji Feng\footnote[1]{fsj@amss.ac.cn} and XiaoSheng Wu\footnote[2]{shengsheng85426@163.com}
\end{center}
\begin{center}
Academy of Mathematics and Systems Science,
Chinese Academy of Sciences, Beijing 100190, P. R. China
\end{center}
\begin{abstract}
Assuming the Riemann Hypothesis, we show that infinitely often consecutive non-trivial zeros of the Riemann
zeta-function differ by at least 2.7327 times the average spacing and infinitely often they differ by at most 0.5154
times the average spacing.
\end{abstract}

\section{Introduction}
Let $\zeta(s)$ denote the Riemann zeta-function. We denote the non-trivial zeros of
$\zeta(s)$ as $\rho=\beta+i\gamma$. Let $\gamma\le\gamma'$ denote consecutive ordinates of the zeros of $\zeta(s)$.
The von Mangoldt formulate (see \cite{Tichmarsh}) gives
\begin{align}
   N(T)=\frac T{2\pi}\log\frac T{2\pi e}+O(\log T),\notag
\end{align}
where $N(T)$ is the number of zeros of $\zeta(s), s=\sigma+it$ in the rectangle $0\le\sigma\le1, 0\le t\le T$. Hence, the average size of $\gamma'-\gamma$ is $2\pi/\log\gamma$. Let
\begin{align}
   \lambda=\lim\sup(\gamma'-\gamma)\frac{\log\gamma}{2\pi}\notag
\end{align}
and
\begin{align}
   \mu=\lim\inf(\gamma'-\gamma)\frac{\log\gamma}{2\pi}\notag,
\end{align}
where $\gamma$ runs over all the ordinates of the zeros of the $\zeta(s)$. In \cite{Montgomery1}, Montgomery suggested that there
exists arbitrarily large and small gaps between consecutive zeros of $\zeta(s)$. That is to say $\mu=0$ and
$\lambda=+\infty$.

Understanding the vertical distribution of the zeros of the zeta-function is very important for a number of reasons. One reason, In particular, in \cite{Conrey3} also \cite{Montgomery3}, it was pointed out that spacing of the zeros of the zeta-function connected with the class number problem for imaginary quadratic fields. Also, the gaps between consecutive zeros of the Riemann zeta-function relate to the zeros of the derivative of the Riemann zeta-function near the critical line, see \cite{Feng1},\cite{Garaev},\cite{Ki},\cite{Soundararajan},\cite{Zhang}.

Unconditionally, in 1946, it's remarked by selberg \cite{Selberg} that $\mu<1<\lambda$. In 2005, by making use of the Wirtinger's inequality and the asymptotic formulae of the fourth mixed moments of the zeta-function and its derivative, R. R. Hall \cite{Hall} proved that $\lambda>2.6306$. This result is even better than what was previously known assuming RH.

Assuming RH,  let $T$ be large and $K=T(\log T)^{-2}$. Further, let
\begin{align}
   h(c)=c-\frac{\textrm{Re}(\sum_{nk\le K}a_k\overline{a_{nk}}g_c(n)\Lambda(n)n^{-1/2})}
   {\sum_{k\le K}\left\vert a_k\right\vert^2}\notag
\end{align}
where
\begin{align}
   g_c(n)=\frac{2\sin(\pi c\frac{\log n}{\log T})}{\pi\log n}\notag
\end{align}
and $\Lambda$ is the von Mangoldt's function. On RH, by an argument using the Guinand-weil explicit formula, Montgomery and Odlyzko \cite{Montgomery2}, in 1981, showed that if $h(c)<1$ for some choice of $c$ and $\{a_n\}$, then $\lambda\ge c$, and if $h(c)>1$ for some choice of $c$ and $\{a_n\}$, then $\mu\le c$. They chose the coefficients
\begin{align}
   a_k=\frac1{k^{\frac1{2}}}f(\frac{\log k}{\log K})\ \ \ \ \textrm{and}\ \ \ \
   a_k=\frac{\lambda(k)}{k^{\frac1{2}}}f(\frac{\log k}{\log K})\notag
\end{align}
where $f$ is a continuous function of bounded variation, and $\lambda(k)$, the Liouville function, equals
$(-1)^{\Omega(k)}$; here, $\Omega(k)$ denotes the total number of primes of $k$. Then they obtain
$\lambda>1.9799$ and $\mu<0.5179$ by optimizing over such functions $f$.

In 1982, by a different method, Mueller \cite{Mueller} showed that $\lambda\geq 1.9$ which is an immediate consequence of a mean value theorem of Gonek \cite{Gonek}.

Montgomery and Odlyzko's results were soon improved by Conrey, Ghosh and Gonek \cite{Conrey1} in 1984. In \cite{Conrey1}, Conrey, Ghosh and Gonek altered  the coefficients to
\begin{align}
   a_k=\frac{d_r(k)}{k^{\frac1{2}}}\ \ \ \ \textrm{and}\ \ \ \
   a_k=\frac{\lambda(k)d_r(k)}{k^{\frac1{2}}}\notag
\end{align}
where $d_r(k)$ is a multiplicative function and for a prime $p$,
\begin{align}
   d_r(p^m)=\frac{\Gamma(m+r)}{\Gamma(r)m!}\notag.
\end{align}
By the choices $r=2.2$ and $r=1.1$ respectively, they improved the values to $\lambda>2.337$ and $\mu<0.5172$.

Recently, H. M. Bui, M. B. Milinovich and N. Ng \cite{Bui2} combined the coefficients of \cite{Montgomery2} and \cite{Conrey1} and got the coefficients
\begin{align}
   a_k=\frac{d_r(k)}{k^{\frac1{2}}}f(\frac{\log k}{\log K})\ \ \ \ \textrm{and}\ \ \ \
   a_k=\frac{\lambda(k)d_r(k)}{k^{\frac1{2}}}f(\frac{\log k}{\log K}).\notag
\end{align}
By optimizing over both $f$ and $r$, they obtain $\lambda>2.69$ and $\mu<0.5155$.

There are still some results about $\lambda$ on the condition of GRH. Conrey, Ghosh and Gonek \cite{Conrey2} combined the ideas of Mueller \cite{Mueller} and \cite{Conrey1} and then proved $\lambda>2.68$, also, by a general mollifier, N. Ng \cite{Ng} improved the result to $\lambda>3$. Still, by a extension of the mollifier of N. Ng, H. M. Bui \cite{Bui1} proved $\lambda>3.0155$.

The works of \cite{Conrey1},\cite{Bui2} are based on the idea of \cite{Montgomery2}. Our work is still based on this, and the results are

{\bf Theorem 1.1}\ \
If the Riemann Hypothesis is true, then  $\lambda>2.7327$ and $\mu<0.5154$.

In order to prove Theorem 1.1, we choose the coefficients
\begin{align}
   a_k=&\frac{d_r(k)}{k^{\frac1{2}}}f_1(\frac{\log K/k}{\log
   K})+\frac{d_r(k)}{k^{\frac1{2}}}\sum_{p_1p_2\mid k}\frac{\mu^2(p_1p_2)\log
   p_1\log p_2}{\log^2K}f_2(\frac{\log K/k}{\log K})\notag
\end{align}
and
\begin{align}
   a_k=&\frac{\lambda(k)d_r(k)}{k^{\frac1{2}}}f_1(\frac{\log K/k}{\log
   K})+\frac{\lambda(k)d_r(k)}{k^{\frac1{2}}}\sum_{p_1p_2\mid k}\frac{\mu^2(p_1p_2)\log
   p_1\log p_2}{\log^2K}f_2(\frac{\log K/k}{\log K}),\notag
\end{align}
where $f_1, f_2$ are some polynomials which will be specified later.

We now give some further insight into the coefficients used here. We note $\sum_{p_1p_2\mid k}\log p_1\log p_2$ approximates the coefficient of $\frac{\zeta'(s)^2}{\zeta(s)}$ and $\lambda(k)\sum_{p_1p_2\mid k}\log p_1\log p_2$ approximates the coefficient of $\frac{\zeta(2s)\zeta'(s)^2}{\zeta(s)^3}$. Hence, The choices of the coefficients used here make $A(t)=\sum_{k\le K}a_kk^{-it}$ behave like a differential polynomial of $\zeta(\frac12+it)^r$ and $\frac{\zeta'}{\zeta}(\frac12+it)$ and a differential polynomial of $\zeta(1+2it)^r/\zeta(\frac12+it)^r$ and $\frac{\zeta'}{\zeta}(\frac12+it)$ respectively.

The coefficients we choose here are motivated by Feng [7]. In [7], the mollifier
\begin{align}
   \psi(s)=&\sum_{k\le y}\frac{\mu(k)}{k^{\frac R{\log T}+s}}(f_1(\frac{\log y/k}{\log
   K})+f_2(\frac{\log y/k}{\log y})\sum_{p_1p_2\mid k}\frac{\log
   p_1\log p_2}{\log^2y}\notag\\
   &+f_3(\frac{\log y/k}{\log y})\sum_{p_1p_2p_3\mid k}\frac{\log
   p_1\log p_2\log p_3}{\log^3y}+\cdots\notag\\
   &+f_I(\frac{\log y/k}{\log y})\sum_{p_1p_2\cdots p_I\mid k}\frac{\log
   p_1\log p_2\cdots\log p_I}{\log^Iy}),\notag
\end{align}
was introduced to improve the result on lower bound of the proportion of zeros of the Riemann zeta-function on the critical line.
\section{Some lemmas}
To prove Theorem 1.1, we need the following lemmas

{\bf Lemma 2.1}\ \
(Mertens Theorem).
\begin{align}
\sum_{p\le y}\frac {\log p}{p}=\log y+O(1).\notag
\end{align}

{\bf Lemma 2.2}\ \
(Levinson \cite{Levinson}).
\begin{align}
\sum_{p\mid j}\frac {\log p}{p}=O(\log \log j).\notag
\end{align}

{\bf Lemma 2.3}\ \
For fixed $r\ge1$,
\begin{align}
   \sum_{k\le x}\frac{d_r(k)^2}{k}=A_r(\log x)^{r^2}+O((\log T)^{r^2-1})\notag
\end{align}
uniformly for $x\le T$.

{\bf Lemma 2.4}\ \
Let $a_i=1,2$ for $1\le i\le m$, and $f$
is a continuous function, $D\ge 1$, then
\begin{align}
   &\int_{1}^{D}\frac {\log ^{a_1-1}x_1}{x_1}dx_1\int_{1}^{\frac
   {D}{x_1}}\frac {\log ^{a_2-1}x_2}{x_2}dx_2\cdots\int_1^{\frac D{x_1x_2\cdots x_m}}\frac{f(x_1x_2\cdots x_mx)}xdx\notag \\
   =&\frac{\prod_{i=1}^m(a_i-1)!}{(\sum_{i=1}^ma_i)!}\int_1^D\frac {f(x)\log
   ^{\sum_{i=1}^ma_i}x}xdx.\notag
\end{align}

Lemma 2.1-Lemma 2.3 are familiar results and Lemma 2.4 is the Lemma 9 of Feng \cite{Feng2}.

\section{Proof of Theorem 1.1}
In this section, we prove the Theorem 1.1. At first, we give a lower bound for $\lambda$ by evaluating
$h(c)$ with the coefficient
\begin{align}
   a_k=&\frac{d_r(k)}{k^{\frac1{2}}}f_1(\frac{\log K/k}{\log
   K})+\frac{d_r(k)}{k^{\frac1{2}}}\sum_{p_1p_2\mid k}\frac{\mu^2(p_1p_2)\log
   p_1\log p_2}{\log^2K}f_2(\frac{\log K/k}{\log K}),\notag
\end{align}
where $r\ge1$ and $f_1$, $f_2$ are polynomials. Employing this coefficient, we have the
denominator in the ratio of sums in the definition of $h(c)$ is
\begin{align}
   \sum_{k\le K}\left\vert a_k\right\vert^2=&\sum_{k\le K}\frac {d_r(k)^2}kf_1(\frac{\log K/k}{\log
   K})^2\notag\\
   &+2\sum_{k\le K}\frac {d_r(k)^2}kf_1(\frac{\log K/k}{\log K})f_2(\frac{\log K/k}{\log K})\sum_{p_1p_2\mid k}\frac{\mu^2(p_1p_2)\log p_1\log p_2}{\log^2K}\notag\\
   &+\sum_{k\le K}\frac {d_r(k)^2}kf_2(\frac{\log K/k}{\log
   K})^2\sum_{p_1p_2\mid k}\frac{\mu^2(p_1p_2)\log p_1\log p_2}{\log^2K}\notag\\
   &\times\sum_{q_1q_2\mid k}\frac{\mu^2(q_1q_2)\log q_1\log q_2}{\log^2K}\notag\\
   =&D_1+D_2+D_3\notag
\end{align}
with obvious meanings. By Abel summation, Lemma 2.3 and recalling that $K=T(\log T) ^{-2}$, we have
\begin{align}
\label{3.1}
   D_1=&A_rr^2\int_1^Kf_1(\frac{\log K/x}{\log K})^2(\log x)^{r^2-1}\frac{dx}{x}+O((\log T)^{r^2-1})\notag\\
   =&A_rr^2(\log K)^{r^2}\int_0^1(1-u)^{r^2-1}f_1(u)^2du+O((\log T)^{r^2-1})\notag\\
   =&A_rr^2(\log T)^{r^2}\int_0^1(1-u)^{r^2-1}f_1(u)^2du+O((\log T)^{r^2-1+\epsilon}),\tag{3.1}
\end{align}
where $\epsilon>0$ is arbitrary and the constant of $O$ is decided by $r$, $\epsilon$ and $f_1$. Substituting $k$ with $p_1p_2k_0$, we find that
\begin{align}
   D_2=&2\sum_{p_1p_2k_0\le K}\frac {\mu^2(p_1p_2)\log p_1\log p_2d_r(p_1p_2k_0)^2}{p_1p_2k_0\log^2K}f_1(\frac{\log K/p_1p_2k_0}{\log K})f_2(\frac{\log K/p_1p_2k_0}{\log K}).\notag
\end{align}
For $d_r(n)=O(n^\epsilon)$ with arbitrary $\epsilon\ge0$,
\begin{align}
   \sum_{i=2}\sum_{p^i\le x}\frac{d_r(p^i)\log p}{p^i}\ll\sum_{p\le x}\log p\sum_{i=2}\frac1{p^{i(1-\epsilon)}}=O(1).\notag
\end{align}
Then we note the terms for which $(k_0,p_1p_2)\neq1$ contribute at most $O((\log T)^{r^2-1})$ in $D_2$. We will also face similar problems in the remainder of the article and not point out any more. By this and Lemma 2.3, we have
\begin{align}
   D_2=&\frac{2r^4}{\log^2K}\sum_{p_1p_2\le K}\frac{\mu^2(p_1p_2)\log p_1\log p_2}{p_1p_2}\sum_{k_0\le K/p_1p_2}\frac{d_r(k_0)^2}{k_0}f_1(\frac{\log K/p_1p_2k_0}{\log K})\notag\\
   &\times f_2(\frac{\log K/p_1p_2k_0}{\log K})+O((\log T)^{r^2-1})\notag\\
   =&\frac{2A_rr^6}{\log^2K}\sum_{p_1p_2\le K}\frac{\mu^2(p_1p_2)\log p_1\log p_2}{p_1p_2}\int_1^{\frac{K}{p_1p_2}}f_1(\frac{\log K/p_1p_2x}{\log K})\notag\\
   &\times f_2(\frac{\log K/p_1p_2x}{\log K})(\log x)^{r^2-1}\frac{dx}{x}+O((\log T)^{r^2-1}).\notag
\end{align}
Let $f$ be a continuous function bounded on $[1,K]$, $i\ge2$ and $a_m\ge1$ are integers for $1\le m\le i$, then
\begin{align}
   &\sum_{p_1p_2\cdots p_i\le K}\frac{\mu^2(p_1p_2\cdots p_i)\log^{a_1}p_1\log^{a_2}p_2\cdots\log^{a_i}p_i}{p_1p_2\cdots p_i}f(p_1p_2\cdots p_i)\notag\\
   =&\sum_{p_1p_2\cdots p_{i-1}\le K}\frac{\mu^2(p_1p_2\cdots p_{i-1})\log^{a_1}p_1\log^{a_2}p_2\cdots\log^{a_{i-1}}p_{i-1}}{p_1p_2\cdots p_{i-1}}\notag\\
   &\times\sum_{{p_i\le K/p_1p_2\cdots p_{i-1}}\atop{(p_i,p_1p_2\cdots p_{i-1})=1}}\frac{\log^{a_i}p_i}{p_i}f(p_1p_2\cdots p_i)\notag\\
   =&\sum_{p_1p_2\cdots p_{i-1}\le K}\frac{\mu^2(p_1p_2\cdots p_{i-1})\log^{a_1}p_1\log^{a_2}p_2\cdots\log^{a_{i-1}}p_{i-1}}{p_1p_2\cdots p_{i-1}}\notag\\
   &\times\sum_{p_i\le K/p_1p_2\cdots p_{i-1}}\frac{\log^{a_i}p_i}{p_i}f(p_1p_2\cdots p_i)\notag\\
   &+\sum_{p_1p_2\cdots p_{i-1}\le K}\frac{\mu^2(p_1p_2\cdots p_{i-1})\log^{a_1}p_1\log^{a_2}p_2\cdots\log^{a_{i-1}}p_{i-1}}{p_1p_2\cdots p_{i-1}}\notag\\
   &\times\sum_{p_i\mid p_1p_2\cdots p_{i-1}}\frac{\log^{a_i}p_i}{p_i}f(p_1p_2\cdots p_i).\notag
\end{align}
By Lemma 2.1 and Lemma 2.2
\begin{align}
   &\sum_{p_1p_2\cdots p_{i-1}\le K}\frac{\mu^2(p_1p_2\cdots p_{i-1})\log^{a_1}p_1\log^{a_2}p_2\cdots\log^{a_{i-1}}p_{i-1}}{p_1p_2\cdots p_{i-1}}\notag\\
   &\times\sum_{p_i\mid p_1p_2\cdots p_{i-1}}\frac{\log^{a_i}p_i}{p_i}f(p_1p_2\cdots p_i)\notag\\
   \ll&\prod_{m=1}^{i-1}\sum_{p\le K}\frac{\log^{a_m}p}{p}\log^{a_i-1}K\log\log K=O((\log K)^{\sum_{m=1}^i a_m-1+\epsilon}).\notag
\end{align}
Then by Lemma 2.1 and Able summation,
\begin{align}
   &\sum_{p_1p_2\cdots p_i\le K}\frac{\mu^2(p_1p_2\cdots p_i)\log^{a_1}p_1\log^{a_2}p_2\cdots\log^{a_i}p_i}{p_1p_2\cdots p_i}f(p_1p_2\cdots p_i)\notag\\
   =&\sum_{p_1p_2\cdots p_{i-1}\le K}\frac{\mu^2(p_1p_2\cdots p_{i-1})\log^{a_1}p_1\log^{a_2}p_2\cdots\log^{a_{i-1}}p_{i-1}}{p_1p_2\cdots p_{i-1}}\notag\\
   &\times\int_1^{\frac{K}{p_1p_2\cdots p_{i-1}}}f(p_1p_2\cdots p_{i-1}x)\log^{a_i-1}x\frac{dx}{x}+O((\log K)^{\sum_{m=1}^i a_m-1+\epsilon}).\notag
\end{align}
We can treat $p_{i-1}$ the same to $p_i$, then by induction, we have
\begin{align}
   &\sum_{p_1p_2\cdots p_i\le K}\frac{\mu^2(p_1p_2\cdots p_i)\log^{a_1}p_1\log^{a_2}p_2\cdots\log^{a_i}p_i}{p_1p_2\cdots p_i}f(p_1p_2\cdots p_i)\notag\\
\label{3.2}
   =&\int_1^K\log^{a_1-1}x_1\frac{dx_1}{x_1}\int_1^{\frac{K}{x_1}}\log^{a_2-1}x_2\frac{dx_2}{x_2}\cdots\int_1^{\frac{K}{x_1x_2\cdots x_{i-1}}}f(x_1x_2\cdots x_i)\log^{a_i-1}x_i\frac{dx_i}{x_i}\notag\\
   &+O((\log K)^{\sum_{m=1}^i a_m-1+\epsilon}),\tag{3.2}
\end{align}
where the constant of $O$ is decided by $f, \epsilon$. Employing this with $i=2$ and
\begin{align}
   f=\frac1{(\log K)^{r^2}}\int_1^{\frac{K}{p_1p_2}}f_1(\frac{\log K/p_1p_2x}{\log K})f_2(\frac{\log K/p_1p_2x}{\log K})(\log x)^{r^2-1}\frac{dx}{x},\notag
\end{align}
we obtain that
\begin{align}
   D_2=&\frac{2A_rr^6}{\log^2K}\int_1^K\frac{dx_1}{x_1}\int_1^{\frac{K}{x_1}}\frac{dx_2}{x_2}\int_1^{\frac{K}{x_1x_2}}f_1(\frac{\log K/x_1x_2x}{\log K})\notag\\
   &\times f_2(\frac{\log K/x_1x_2x}{\log K})(\log x)^{r^2-1}\frac{dx}{x}+O((\log T)^{r^2-1+\epsilon}).\notag
\end{align}
By Lemma 2.4, we have
\begin{align}
   D_2=&\frac{2A_rr^6}{\log^2K}\int_1^K\log x_1\frac{dx_1}{x_1}\int_1^{\frac{K}{x_1}}f_1(\frac{\log K/x_1x}{\log K})f_2(\frac{\log K/x_1x}{\log K})(\log x)^{r^2-1}\frac{dx}{x}\notag\\
   &+O((\log T)^{r^2-1+\epsilon}).\notag
\end{align}
By variable changes $u=1-\frac{\log x_1}{\log K}$, $v=1-\frac{\log x_1x}{\log K}$,
\begin{align}
\label{3.3}
   D_2=&2A_rr^6(\log K)^{r^2}\int_0^1(1-u)\int_0^u(u-v)^{r^2-1}f_1(v)f_2(v)dvdu+O((\log T)^{r^2-1+\epsilon})\notag\\
   =&2A_rr^6(\log T)^{r^2}\int_0^1(1-u)\int_0^u(u-v)^{r^2-1}f_1(v)f_2(v)dvdu+O((\log T)^{r^2-1+\epsilon}),\tag{3.3}
\end{align}
where the constant of $O$ is decided by $r$, $\epsilon$ and $f_1$, $f_2$. It's easy to find (see Lemma 8 of \cite{Feng2})
\begin{align}
\label{3.31}
   &\sum_{p_1p_2\mid k}\mu^2(p_1p_2)\log p_1\log p_2\sum_{q_1q_2\mid k}\mu^2(q_1q_2)\log q_1\log q_2\notag\\
   =&\sum_{p_1p_2p_3p_4\mid k}\mu^2(p_1p_2p_3p_4)\log p_1\log p_2\log p_3\log p_4\notag\\
   &+4\sum_{p_1p_2p_3\mid k}\mu^2(p_1p_2p_3)\log^2p_1\log p_2\log p_3+2\sum_{p_1p_2\mid k}\mu^2(p_1p_2)\log^2p_1\log^2p_2.\tag{3.4}
\end{align}
By this, we have
\begin{align}
   D_3=&\sum_{k\le K}\frac {d_r(k)^2}kf_2(\frac{\log K/k}{\log
   K})^2\sum_{p_1p_2p_3p_4\mid k}\frac{\mu^2(p_1p_2p_3p_4)\log p_1\log p_2\log p_3\log p_4}{\log^4K}\notag\\
   &+4\sum_{k\le K}\frac {d_r(k)^2}kf_2(\frac{\log K/k}{\log
   K})^2\sum_{p_1p_2p_3\mid k}\frac{\mu^2(p_1p_2p_3)\log^2p_1\log p_2\log p_3}{\log^4K}\notag\\
   &+2\sum_{k\le K}\frac {d_r(k)^2}kf_2(\frac{\log K/k}{\log
   K})^2\sum_{p_1p_2\mid k}\frac{\mu^2(p_1p_2)\log^2p_1\log^2p_2}{\log^4K}\notag\\
   =&D_{31}+D_{32}+D_{33}\notag
\end{align}
with obvious meanings. Similarly to $D_2$, interchanging the summations, by (\ref{3.2}), Lemma 2.3 and Lemma 2.4, we have
\begin{align}
   D_{31}=&\frac{r^8}{\log^4K}\sum_{p_1p_2p_3p_4\le K}\frac{\mu^2(p_1p_2p_3p_4)\log p_1\log p_2\log p_3\log p_4}{p_1p_2p_3p_4}\notag\\
   &\times\sum_{k\le K/p_1p_2p_3p_4}\frac{d_r(k)^2}{k}f_2(\frac{\log K/p_1p_2p_3p_4k}{\log K})^2+O((\log T)^{r^2-1})\notag\\
   =&\frac{A_rr^{10}}{\log^4K}\sum_{p_1p_2p_3p_4\le K}\frac{\mu^2(p_1p_2p_3p_4)\log p_1\log p_2\log p_3\log p_4}{p_1p_2p_3p_4}\notag\\
   &\times\int_1^{\frac{K}{p_1p_2p_3p_4}}f_2(\frac{\log K/p_1p_2p_3p_4x}{\log K})^2(\log x)^{r^2-1}\frac{dx}{x}+O((\log T)^{r^2-1})\notag\\
   =&\frac{A_rr^{10}}{6\log^4K}\int_1^K\log^3x_1\frac{dx_1}{x_1}\int_1^{\frac{K}{x_1}}f_2(\frac{\log K/x_1x}{\log K})^2(\log x)^{r^2-1}\frac{dx}{x}+O((\log T)^{r^2-1+\epsilon})\notag\\
   =&\frac1{6}r^{10}A_r(\log T)^{r^2}\int_0^1(1-u)^3\int_0^u(u-v)^{r^2-1}f_2(v)^2dvdu+O((\log T)^{r^2-1+\epsilon}),\notag
\end{align}
with variable changes $u=1-\frac{\log x_1}{\log K}, v=1-\frac{\log x_1x}{\log K}$. Similarly,
\begin{align}
   D_{32}=&\frac{4r^6}{\log^4K}\sum_{p_1p_2p_3\le K}\frac{\mu^2(p_1p_2p_3)\log^2p_1\log p_2\log p_3}{p_1p_2p_3}\notag\\
   &\times\sum_{k\le K/p_1p_2p_3}\frac{d_r(k)^2}{k} f_2(\frac{\log K/p_1p_2p_3k}{\log K})^2+O((\log T)^{r^2-1})\notag\\
   =&\frac{4A_rr^{8}}{\log^4K}\sum_{p_1p_2p_3\le K}\frac{\mu^2(p_1p_2p_3)\log^2p_1\log p_2\log p_3}{p_1p_2p_3}\notag\\
   &\times\int_1^{\frac{K}{p_1p_2p_3}}f_2(\frac{\log K/p_1p_2p_3x}{\log K})^2(\log x)^{r^2-1}\frac{dx}{x}+O((\log T)^{r^2-1})\notag\\
   =&\frac{2A_rr^{8}}{3\log^4K}\int_1^K\log^3x_1\frac{dx_1}{x_1}\int_1^{\frac{K}{x_1}}f_2(\frac{\log K/x_1x}{\log K})^2(\log x)^{r^2-1}\frac{dx}{x}+O((\log T)^{r^2-1+\epsilon})\notag\\
   =&\frac2{3}r^{8}A_r(\log T)^{r^2}\int_0^1(1-u)^3\int_0^u(u-v)^{r^2-1}f_2(v)^2dvdu+O((\log T)^{r^2-1+\epsilon}),\notag
\end{align}
and
\begin{align}
   D_{33}=&\frac{2r^4}{\log^4K}\sum_{p_1p_2\le K}\frac{\mu^2(p_1p_2)\log^2p_1\log^2p_2}{p_1p_2}\sum_{k\le K/p_1p_2}\frac{d_r(k)^2}{k}f_2(\frac{\log K/p_1p_2k}{\log K})^2\notag\\
   &+O((\log T)^{r^2-1})\notag\\
   =&\frac{2A_rr^{6}}{\log^4K}\sum_{p_1p_2\le K}\frac{\mu^2(p_1p_2)\log^2p_1\log^2p_2}{p_1p_2}\int_1^{\frac{K}{p_1p_2}}f_2(\frac{\log K/p_1p_2x}{\log K})^2(\log x)^{r^2-1}\frac{dx}{x}\notag\\
   &+O((\log T)^{r^2-1})\notag\\
   =&\frac{A_rr^{6}}{3\log^4K}\int_1^K\log^3x_1\frac{dx_1}{x_1}\int_1^{\frac{K}{x_1}}f_2(\frac{\log K/x_1x}{\log K})^2(\log x)^{r^2-1}\frac{dx}{x}+O((\log T)^{r^2-1+\epsilon})\notag\\
   =&\frac1{3}r^{6}A_r(\log T)^{r^2}\int_0^1(1-u)^3\int_0^u(u-v)^{r^2-1}f_2(v)^2dvdu+O((\log T)^{r^2-1+\epsilon}).\notag
\end{align}
Hence,
\begin{align}
\label{3.4}
   D_3=&A_r(\log T)^{r^2}(\frac1{6}r^{10}+\frac2{3}r^8+\frac1{3}r^6)\int_0^1(1-u)^3\int_0^u(u-v)^{r^2-1}f_2(v)^2dvdu\notag\\
   &+O((\log T)^{r^2-1+\epsilon}),\tag{3.5}
\end{align}
where the constant of $O$ is decided by $r$, $\epsilon$ and $f_2$.

We now evaluate the numerator in the ratio of sums in the definition
of $h(c)$. If we let
\begin{align}
   N(c)=\sum_{nk\le K}a_ka_{nk}g_c(n)\Lambda(n)n^{-1/2},\notag
\end{align}
then a straightforward argument shows that
\begin{align}
   N(c)=&\frac2{\pi}\sum_{nk\le
   K}\frac{d_r(k)d_r(nk)\Lambda(n)}{kn\log n}\sin(\pi c\frac{\log
   n}{\log T})(f_1(\frac{\log K/k}{\log K})f_1(\frac{\log K/nk}{\log
   K})\notag\\
   &+f_1(\frac{\log K/nk}{\log K})f_2(\frac{\log K/k}{\log
   K})\sum_{p_1p_2\mid k}\frac{\mu^2(p_1p_2)\log p_1\log p_2}{\log^2K}\notag\\
   &+f_1(\frac{\log K/k}{\log K})f_2(\frac{\log K/nk}{\log
   K})\sum_{p_1p_2\mid nk}\frac{\mu^2(p_1p_2)\log p_1\log p_2}{\log^2K}\notag\\
   &+f_2(\frac{\log K/k}{\log
   K})f_2(\frac{\log K/nk}{\log K})\sum_{p_1p_2\mid k}\frac{\mu^2(p_1p_2)\log p_1\log p_2}{\log^2K}\notag\\
   &\times\sum_{q_1q_2\mid nk}\frac{\mu^2(q_1q_2)\log q_1\log q_2}{\log^2K})\notag\\
   =&N_1+N_2+N_3+N_4\notag
\end{align}
with the obvious meaning. By the familiar distribution of $\Lambda(n)$, we have
\begin{align}
   N_1=&\frac2{\pi}\sum_{nk\le
   K}\frac{d_r(k)d_r(nk)\Lambda(n)}{kn\log n}\sin(\pi c\frac{\log
   n}{\log T})f_1(\frac{\log K/k}{\log K})f_1(\frac{\log K/nk}{\log
   K})\notag\\
   =&\frac2{\pi}\sum_{pk\le K}\frac{d_r(k)d_r(pk)}{kp}\sin(\pi c\frac{\log
   p}{\log T})f_1(\frac{\log K/k}{\log K})f_1(\frac{\log K/pk}{\log
   K})+O((\log T)^{r^2-1})\notag\\
   =&\frac{2r}{\pi}\sum_{p\le K}\frac{\sin(\pi c\frac{\log
   p}{\log T})}{p}\sum_{k\le K/p}\frac{d_r(k)^2}{k}f_1(\frac{\log K/k}{\log K})f_1(\frac{\log K/pk}{\log
   K})+O((\log T)^{r^2-1}).\notag
\end{align}
By Lemma 2.3 and Abel summation, we have
\begin{align}
   N_1=&\frac{2A_rr^3}{\pi}\sum_{p\le K}\frac{\sin(\pi c\frac{\log p}{\log T})}{p}
   \int_1^{\frac K{p}}f_1(\frac{\log K/x}{\log K})f_1(\frac{\log K/px}{\log K})(\log x)^{r^2-1}\frac{dx}{x}\notag\\
   &+O((\log T)^{r^2-1}).\notag
\end{align}
From Lemma 2.1 and Abel summation,
\begin{align}
   N_1=&\frac{2A_rr^3}{\pi}\int_1^K\frac{\sin{(\pi c\frac{\log x_1}{\log T})}}{x_1\log x_1}dx_1\int_1^{\frac K{x_1}}
   f_1(\frac{\log K/x}{\log K})f_1(\frac{\log K/xx_1}{\log K})(\log x)^{r^2-1}\frac{dx}{x}\notag\\
   &+O((\log T)^{r^2-1}).\notag
\end{align}
Interchanging the order of integration, it follows that
\begin{align}
   N_1=&\frac{2A_rr^3}{\pi}\int_1^Kf_1(\frac{\log K/x_1}{\log K})(\log x_1)^{r^2-1}\frac{dx_1}{x_1}\int_1^{\frac K{x_1}}
   \frac{\sin{(\pi c\frac{\log x}{\log T})}}{\log x}f_1(\frac{\log K/xx_1}{\log K})\frac{dx}{x}\notag\\
   &+O((\log T)^{r^2-1}).\notag
\end{align}
Let $u=1-\frac{\log x_1}{\log K}$, $v=\frac{\log x}{\log K}$, we have
\begin{align}
\label{3.5}
   N_1=&\frac{2A_rr^3}{\pi}(\log T)^{r^2}\int_0^1(1-u)^{r^2-1}f_1(u)\int_0^u\frac{\sin(\pi
   cv\frac{\log K}{\log T})}{v}f_1(u-v)dvdu\notag\\
   &+O((\log T)^{r^2-1+\epsilon})\notag\\
   =&\frac{2A_rr^3}{\pi}(\log T)^{r^2}\int_0^1(1-u)^{r^2-1}f_1(u)\int_0^u\frac{\sin(\pi
   cv)}{v}f_1(u-v)dvdu\notag\\
   &+O((\log T)^{r^2-1+\epsilon}),\tag{3.6}
\end{align}
where the constant of $O$ is decided by $r$, $\epsilon$ and $f_1$. Similarly,
\begin{align}
   N_2=&\frac2{\pi}\sum_{pk\le K}\sin(\pi c\frac{\log
   p}{\log T})\frac{d_r(k)d_r(kp)}{kp}f_1(\frac{\log K/pk}{\log K})f_2(\frac{\log K/k}{\log
   K})\notag\\
   &\times\sum_{p_1p_2\mid k}\frac{\mu^2(p_1p_2)\log p_1\log p_2}{\log^2K}+O((\log T)^{r^2-1})\notag\\
   =&\frac{2r^5}{\pi\log^2K}\sum_{p_1p_2\le K}\frac{\mu^2(p_1p_2)\log p_1\log p_2}{p_1p_2}\sum_{pk\le K/p_1p_2}\frac{\sin(\pi c\frac{\log
   p}{\log T})d_r(k)^2}{pk}\notag\\
   &\times f_1(\frac{\log K/pp_1p_2k}{\log K})f_2(\frac{\log K/p_1p_2k}{\log
   K})+O((\log T)^{r^2-1})\notag
\end{align}
From the calculation of $N_1$, we know the inner sum in the main term of the last expression of $N_2$ is
\begin{align}
   &A_rr^2\int_1^{\frac K{p_1p_2}}f_2(\frac{\log K/p_1p_2x_1}{\log K})(\log x_1)^{r^2-1}\frac{dx_1}{x_1}\notag\\
   &\times\int_1^{\frac K{p_1p_2x_1}}
   \frac{\sin{(\pi c\frac{\log x}{\log T})}}{\log x}f_1(\frac{\log K/p_1p_2xx_1}{\log K})\frac{dx}{x}+O((\log T)^{r^2-1}).\notag
\end{align}
Employing this and by (\ref{3.2}), an argument similar to $D_2$ shows that
\begin{align}
   N_2=&\frac{2A_rr^7}{\pi(\log K)^2}\int_1^K\log x_1\frac{dx_1}{x_1}\int_1^{\frac K{x_1}}f_2(\frac{\log K/x_1x_2}{\log K})(\log x_2)^{r^2-1}\frac{dx_2}{x_2}\notag\\
   &\times\int_1^{\frac K{x_1x_2}}\frac{\sin{(\pi c\frac{\log x}{\log T})}}{\log x}f_1(\frac{\log K/xx_1x_2}{\log K})\frac{dx}{x}+O((\log T)^{r^2-1}).\notag
\end{align}
Then, by variable changes $u=1-\frac{\log x_1}{\log K}, v=1-\frac{\log x_1x_2}{\log K}, w=\frac{\log x}{\log K}$,
\begin{align}
\label{3.6}
   N_2=&\frac{2A_rr^7}{\pi}(\log T)^{r^2}\int_0^1(1-u)\int_0^u(u-v)^{r^2-1}f_2(v)\int_0^v\frac{\sin(\pi cw)}{w}f_1(v-w)dwdvdu\notag\\
   &+O((\log T)^{r^2-1+\epsilon}),\tag{3.7}
\end{align}
where the constant of $O$ is decided by $r$, $\epsilon$ and $f_1$, $f_2$. Still,
\begin{align}
   N_3=&\frac2{\pi\log^2K}\sum_{pk\le K}\sin(\pi c\frac{\log
   p}{\log T})\frac{d_r(k)d_r(kp)}{kp}f_1(\frac{\log K/k}{\log K})f_2(\frac{\log
   K/pk}{\log K})\notag\\
   &\times\sum_{p_1p_2\mid pk}\mu^2(p_1p_2)\log p_1\log p_2+O((\log T)^{r^2-1})\notag.
\end{align}
A simple calculation shows that
\begin{align}
\label{3.7}
   \sum_{p_1p_2\mid pk}\mu^2(p_1p_2)\log p_1\log p_2=\sum_{p_1p_2\mid k}\mu^2(p_1p_2)\log p_1\log p_2+2\log p\sum_{p_1\mid k}\log p_1,\tag{3.8}
\end{align}
for $(p,k)=1$. For those items with $(p,k)\neq1$ in $N_3$ contribute $O((\log T)^{r^2-1})$ at most, therefore, we can employ this in the main term of the expression of $N_3$ and denote that $N_3=N_{31}+N_{32}+O((\log T)^{r^2-1})$, where
\begin{align}
   N_{31}=&\frac{2r^5}{\pi\log^2K}\sum_{p_1p_2\le K}\frac{\mu^2(p_1p_2)\log p_1\log p_2}{p_1p_2}\sum_{pk\le K/p_1p_2}\frac{\sin(\pi c\frac{\log
   p}{\log T})d_r(k)^2}{pk}\notag\\
   &\times f_2(\frac{\log K/pp_1p_2k}{\log K})f_1(\frac{\log K/p_1p_2k}{\log
   K})+O((\log T)^{r^2-1})\notag\\
   =&\frac{2A_rr^7}{\pi(\log K)^2}\int_1^K\log x_1\frac{dx_1}{x_1}\int_1^{\frac K{x_1}}f_1(\frac{\log K/x_1x_2}{\log K})(\log x_2)^{r^2-1}\frac{dx_2}{x_2}\notag\\
   &\times\int_1^{\frac K{x_1x_2}}\frac{\sin{(\pi c\frac{\log x}{\log T})}}{\log x}f_2(\frac{\log K/xx_1x_2}{\log K})\frac{dx}{x}+O((\log T)^{r^2-1+\epsilon})\notag\\
   =&\frac{2A_rr^7}{\pi}(\log T)^{r^2}\int_0^1(1-u)\int_0^u(u-v)^{r^2-1}f_1(v)\notag\\
   &\times\int_0^v\frac{\sin(\pi cw)}{w}f_2(v-w)dwdvdu+O((\log T)^{r^2-1+\epsilon}),\notag
\end{align}
and
\begin{align}
   N_{32}=&\frac{4r^3}{\pi\log^2K}\sum_{p_1\le K}\frac{\log p_1}{p_1}\sum_{pk\le K/p_1}\frac{\sin(\pi c\frac{\log
   p}{\log T})d_r(k)^2\log p}{pk}f_2(\frac{\log K/pp_1k}{\log K})\notag\\
   &\times f_1(\frac{\log K/p_1k}{\log
   K})+O((\log T)^{r^2-1})\notag\\
   =&\frac{4A_rr^5}{\pi(\log K)^2}\int_1^K\frac{dx_1}{x_1}\int_1^{\frac K{x_1}}f_1(\frac{\log K/x_1x_2}{\log K})(\log x_2)^{r^2-1}\frac{dx_2}{x_2}\int_1^{\frac K{x_1x_2}}\sin{(\pi c\frac{\log x}{\log T})}\notag\\
   &\ \ \ \times f_2(\frac{\log K/xx_1x_2}{\log K})\frac{dx}{x}+O((\log T)^{r^2-1+\epsilon})\notag\\
   =&\frac{4A_rr^5}{\pi}(\log T)^{r^2}\int_0^1\int_0^u(u-v)^{r^2-1}f_1(v)\int_0^v\sin(\pi cw)f_2(v-w)dwdvdu\notag\\
   &\ \ \ +O((\log T)^{r^2-1+\epsilon}),\notag
\end{align}
So, we have
\begin{align}
   N_3=&\frac{2A_rr^7}{\pi}(\log T)^{r^2}
   \int_0^1(1-u)\int_0^u(u-v)^{r^2-1}f_1(v)\int_0^v\frac{\sin(\pi cw)}{w}f_2(v-w)dwdvdu\notag\\
\label{3.8}
   &+\frac{4A_rr^5}{\pi}(\log T)^{r^2}\int_0^1\int_0^u(u-v)^{r^2-1}f_1(v)\int_0^v\sin(\pi cw)f_2(v-w)dwdvdu\notag\\
   &+O((\log T)^{r^2-1+\epsilon}),\tag{3.9}
\end{align}
where the constant of $O$ is decided by $r$, $\epsilon$ and $f_1$, $f_2$. Finally,
\begin{align}
  N_4=&\frac2{\pi\log^4K}\sum_{pk\le K}\sin(\pi c\frac{\log
   p}{\log T})\frac{d_r(k)d_r(kp)}{kp}f_2(\frac{\log K/k}{\log K})f_2(\frac{\log K/pk}{\log K})\notag\\
   &\times\sum_{p_1p_2\mid k}\mu^2(p_1p_2)\log p_1\log p_2\sum_{q_1q_2\mid pk}\mu^2(q_1q_2)\log q_1\log q_2
   +O((\log T)^{r^2-1}).\notag
\end{align}
By (\ref{3.31}) and (\ref{3.7}), we have
\begin{align}
   &\sum_{p_1p_2\mid k}\mu^2(p_1p_2)\log p_1\log p_2\sum_{q_1q_2\mid pk}\mu^2(q_1q_2)\log q_1\log q_2\notag\\
   =&\sum_{p_1p_2p_3p_4\mid k}\mu^2(p_1p_2p_3p_4)\log p_1\log p_2\log p_3\log p_4\notag\\
   &+4\sum_{p_1p_2p_3\mid k}\mu^2(p_1p_2p_3)\log^2p_1\log p_2\log p_3+2\sum_{p_1p_2\mid k}\mu^2(p_1p_2)\log^2p_1\log^2p_2\notag\\
   &+2\log p\sum_{p_1p_2p_3\mid k}\mu^2(p_1p_2p_3)\log p_1\log p_2\log p_3+4\log p\sum_{p_1p_2\mid k}\mu^2(p_1p_2)\log^2p_1\log p_2,\notag
\end{align}
for $(k,p)=1$. Hence, an argument as before shows that
\begin{align}
   N_4=&\frac{A_r}{\pi}(\frac1{3}r^{11}+\frac4{3}r^9+\frac2{3}r^7)(\log T)^{r^2}\int_0^1(1-u)^3\int_0^u(u-v)^{r^2-1}f_2(v)\notag\\
   &\times\int_0^v\frac{\sin(\pi cw)}{w}f_2(v-w)dwdvdu\notag\\
   &+\frac{A_r}{\pi}(2r^{9}+4r^7)(\log T)^{r^2}\int_0^1(1-u)^2\int_0^u(u-v)^{r^2-1}f_2(v)\notag\\
   &\times\int_0^v\sin(\pi cw)f_2(v-w)dwdvdu\notag\\
\label{3.9}
   &+O((\log T)^{r^2-1+\epsilon}),\tag{3.10}
\end{align}
where the constant of $O$ is decided by $r$, $\epsilon$ and $f_2$. Consequently, we find that
\begin{align}
   h(c)=c-\frac{N_1+N_2+N_3+N_4}{D_1+D_2+D_3},\notag
\end{align}
where $D_1-D_3$, $N_1-N_4$ are given by (\ref{3.1}), (\ref{3.3}), (\ref{3.4})-(\ref{3.6}), (\ref{3.8}), (\ref{3.9}). Choosing $r=2.6$ and $f_1(x)=-3.54-42.94x+88.05x^2-34.33x^3$
, $f_2(x)=4.56+63.02x+42.72x^2+34.45x^3$, we obtain (by a numerical calculation) that $h(2.7327)<1$
when $T$ is sufficiently large. This provides the lower bound for $\lambda$ in Theorem 1.1.

Since $\lambda(n)^2=1$ and $\lambda(np)=-\lambda(n)$ for every $n\in N$ and every prime $p$, we can evaluate $h(c)$
with the coefficient
\begin{align}
   a_k=&\frac{\lambda(k)d_r(k)}{k^{\frac1{2}}}f_1(\frac{\log K/k}{\log
   K})+\frac{\lambda(k)d_r(k)}{k^{\frac1{2}}}\sum_{p_1p_2\mid k}\frac{\mu^2(p_1p_2)\log
   p_1\log p_2}{\log^2K}f_2(\frac{\log K/k}{\log K})\notag
\end{align}
as before. Here, as above, $f_1, f_2$ are polynomials. With this
choice of coefficient, we obtain that
\begin{align}
   h(c)=c+\frac{N_1+N_2+N_3+N_4}{D_1+D_2+D_3},\notag
\end{align}
where $D_1-D_3$, $N_1-N_4$ are given by (\ref{3.1}), (\ref{3.3}), (\ref{3.4})-(\ref{3.6}), (\ref{3.8}), (\ref{3.9}). Choosing $r=1.18$ and $f_1(x)=1.25+0.95x+2.07x^2-2.21x^3$
, $f_2(x)=0.7+1.92x$, we obtain (by a numerical calculation) that $h(0.5154)>1$
when $T$ is sufficiently large. This provides the lower bound for $\mu$ in Theorem 1.1.

It's worth to remark that we may generalize the coefficient to
\begin{align}
   a_k=&\frac{d_r(k)}{k^{\frac1{2}}}(f_1(\frac{\log K/k}{\log
   K})+f_2(\frac{\log K/k}{\log K})\sum_{p_1p_2\mid k}\frac{\log
   p_1\log p_2}{\log^2K}\notag\\
   &+f_3(\frac{\log K/k}{\log K})\sum_{p_1p_2p_3\mid k}\frac{\log
   p_1\log p_2\log p_3}{\log^3K}+\cdots\notag\\
   &+f_I(\frac{\log K/k}{\log K})\sum_{p_1p_2\cdots p_I\mid k}\frac{\log
   p_1\log p_2\cdots\log p_I}{\log^IK}),\notag
\end{align}
for any integer $I\ge2$. There is no problem to calculate $h(c)$ with this coefficient. However, the numerical calculation doesn't seem
like there is much to gain by increasing any more item of the coefficient we choose.

\end{document}